\documentclass[twoside]{article}
\def\Vox{$^{\fbox {\/}}$}
\def\prom{
\newtheorem{prop}{Proposition}[section]

\newtheorem{lem}[prop]{Lemma}
\newtheorem{exa}[prop]{Example}

\newtheorem{rem}[prop]{Remark}

\newtheorem{tef}[prop]{Definition}
}
\def\ieacosdo{
\author {\ \\N.\ C.\ A.\ da Costa\\F.\ A.\ Doria\\\
\\Institute for Advanced
Studies, University of S\~ao Paulo.\\Av.\ Prof.\
Luciano Gualberto, trav.\ J, 374.\\05655--010 S\~ao
Paulo SP Brazil.\\\ \\{\sc ncacosta@usp.br}\\{\sc
doria1000@yahoo.com.br}\\{\sc doria@lncc.br}}}
\prom

\begin {document}
\pagestyle {myheadings}
\title {On a total function which overtakes all
total recursive functions.\thanks {Partially
supported by FAPESP and CNPq. Alternative address
for F.\ A.\ Doria: Research Center for Mathematical
Theories of Communication and Program IDEA, School
of Communications, Federal University at Rio de
Janeiro, Av.\ Pasteur, 250. 22295--900 Rio RJ
Brazil.}}
\ieacosdo
\maketitle

\begin {abstract}\noindent We discuss here a function
that is part of the folklore around the $P=?NP$
problem. This function $g^*$ is defined over all
time--polynomial Turing machines; we trivially change
it into a $g$ which is defined over all Turing
machines, and show that $g$ overtakes all total
recursive functions. We then show that $g$ is
dominated by a generalized Busy Beaver function.
\end {abstract}

\newpage
\thispagestyle {plain}

\markboth {da Costa, Doria}{Busy Beaver}
\section {Introduction}
\markboth {da Costa, Doria}{Busy Beaver}

Mathematics is sometimes full of unexpected turns,
like a good mystery novel. The Busy Beaver game
looks na\"\i ve enough, and yet it leads us to an
incredibly fast--growing, uncomputable function
\cite {Rado}. Also, who would expect that questions
about finite functions would depend \cite {Fried}
on the existence of large cardinals? 

Sometimes mathematics behaves in a rather
predictable, reasonable way. When one goes from real
numbers to complex numbers to quaternions and then
to the Cayley numbers, each step in that succession
is marked (as it should be) by a slight modification
of previously known properties. Yet in many cases
the generalization of properties leads us into the
unexpected: a contemporary example can be found in
the mathematics of gauge fields. Given an Abelian
gauge field, that field is always (locally, at
least) derived from a gauge potential (a connection
form) which is unique modulo gauge transformations.
However when we deal with non--Abelian gauge fields,
we see that the same field can be derived from
physically different gauge potentials, that is,
potentials that aren't related by a gauge
transformation, not even locally \cite {Do1,Do2,San}.
Just by looking at the Abelian case we would never
expect such a novel phenomenon to creep up in the
non--Abelian case. Moreover this precise phenomenon
is quite surprisingly related to a kind of
transformation first used by Einstein in his theory
of the asymmetric field
\cite {Do1,Eins}.  

Can we say that we are dealing here with
singularities in a kind of space of mathematical
theories? We have in fact recently proposed this
picture: we described a coding of formal systems by
dynamical systems that bifurcate: for instance, if a
given (undecidable) sentence holds, the system
bifurcates; if not, it doesn't. Thus undecidable
sentences are coded by singularities in those
systems \cite {CosDo4}. 

(However we are not sure whether the formal picture
that opposes singular points to regular points is an
apt characterization of our intuitions about the
unexpected in mathematics; the picture may be still
more complex.) 

\subsection*{The $P=?NP$ question} 

Let's consider a final example, which will lead us
into this paper's main topic. That example stems from
the $P=?NP$ question. 

The $P=?NP$ question starts from a very simple,
rather obvious question that we can formulate as a
short tale:
\begin {quote} Mrs.\ H.\ is a gentle and able lady
who has long been the secretary of a large
university department. Every semester Mrs.\ H.\ is
confronted with the following problem: there are
courses to be taught, professors to be distributed
among different classes of students, large and small
classes, and a shortage of classrooms. She fixes a
minimum acceptable level of overlap among classes
and students and sets down in a tentative way to get
the best possible schedule given that minimum
desired overlap. It's a tiresome task, and in most
cases, when there are many new professors or when the
dean changes the classroom allocation system, Mrs.\
H.\ quite frequently has to check nearly all 
conceivable schedulings before she is able to reach a
conclusion. In despair she asks a professor whom she
knows has a degree in math: {\em ``tell me, can't
you find in your math a fast way of scheduling our
classes with a minimum level of overlap among
them?''}
\end {quote}

Mrs.\ H.\ unknowingly asks about the $P=?NP$
question. She is able to understand its basic
contours: {\em there are questions for which it is
easy to check whether a given arrangement of data
fits in as a solution; however for the general case
there are no known shortcuts in order to reach a
solution.} 

If you can guess the answer, it is easy to verify
it. If you have to work it out, you'll most likely
have to sweat your shirt until you reach an adequate
answer. Most likely: for one always wonders whether
there is some general way to get a quick solution.
Who knows? 

Are there general shortcuts available? This is the
whole point---in a nutshell; and this is the main
query in the $P=?NP$ problem. 

$P=NP$ is, roughly, ``there is a time--polynomial
Turing machine that inputs each instance of a
problem in  a $NP$ class and that correctly guesses
(outputs) a solution for that instance.'' The
negation $P<NP$ is, ``for each polynomial machine of
G\"odel number $m$ there is an instance $x$ such
that machine $m$ guesses wrongly about $x$.'' 

The counterexample function $f(m)$ is the function
that enumerates each first instance $x$ where a
polynomial algorithm fails; it is recursive on the
set of all polynomial machines and $P<NP$ holds if
and only if $f$ is total. 

So, one of the possible approaches to this major
problem is to focus on the counterexample function
$f$. How are we to understand its properties?
Perhaps by analogy: we use a tamer function $g$ that
more or less looks like $f$, and try to infer the
properties of $f$ from those of $g$. The chief
property we want to consider is whether $f$ can be
proved to be total in Peano Arithmetic (PA), for in
that case: 
\begin {quote}
If the counterexample function $f$ is bounded by a
prescribed total recursive function---for Peano
Arithmetic, the fast--growing function ${\sf
F}_{\epsilon_0}$---then $P<NP$ will be provable in
Peano Arithmetic. 
\end {quote}

For an adequately conceived $g$, it should be easy
to check the desired property for $g$, and it should
now be just a small step from $g$ to $f$. 

(For ${\sf F}_{\epsilon_0}$ see \cite {Spen} and
references therein; for the counterexample function
see \cite {Do}.) 

Then our main query is: is $g$ bounded by such a
function ${\sf F}_{\epsilon_0}$? If not, where do the
properties of $g$ will lead us concerning $f$\,?  

\subsection*{Functions ${\sf g}_0$ and $g^*$}

\begin {rem}\rm 
Let's motivate $g^*$. Consider the
following: 
\begin {equation}\label {1}
{\sf g}_0(m) = \mu _x (x^m < 2^x)
\end {equation} for
non--negative integers $m$. This function is
trivially recursive and total, and bounded by the
exponential $2^x$ itself. 

Now suppose that we have defined some (necessarily
nonrecursive) enumeration for all time--polynomial
Turing machines (the poly machines). This
enumeration is noted,  ${\sf P}_m, m\in\omega$; we
note ${\cal P}$ the set of all polymachines. The
following function:
\begin {equation}\label {2}
g^*(m) = \mu_x [{\sf P}_m (x) < 2^x]
\end {equation}
looks very similar to (\ref {1}). In fact we have
only changed the monomial $x^m$ for the output
${\sf P}_m (x)$ of a poly machine, and that output
is necessarily bounded in its length by some
polynomial, that is, there is a positive integer $k$
so that
$$|{\sf P}_m(x)|<|x|^k.$$

Therefore $g^*$ is also total. 

$g^*$ is usually presented as a kind of analogue of
${\sf g}_0$ in eq.\ (\ref {1}), and it is argued
that, since ${\sf g}_0$ is bounded by an exponential,
the same should be true of $g^*$ in eq.\ (\ref {2}).
Yes, a moment's thought shows that this should be
the case\ldots 

Wait a moment---is it really so?

Let's ponder it: the relation between exponent $k$
and machine index $m$ is possibly quite
complicated---and index $m$ roams over a
nonrecursive set. 

The analogy starts to break down here. 
 
The analogy breaks down right at the beginning. Can
we pass over those differences and still have a
rather ``tame'' behavior for $g^*$ as the one
exhibited by ${\sf g}_0$? \Vox
\end {rem}

We show in this paper how wildly different is the
behavior of $g^*$ when compared to ${\sf g}_0$,
despite the fact that the two functions are
apparently very similar. The whole point, as we will
see, is the relation between machine index $m$ and
the exponent $k$ for the bounding polynomial. 

\subsection*{Analogies between $f$ and $g^*$}

The main analogy between the counterexample function
$f$ and $g^*$ is: both are defined over the set of
all poly machines ${\cal P}$ and their values are
given by the application of the $\mu$ operator to
a recursive predicate on ${\cal P}$ that depends on
the index $m$ of a poly machine ${\sf P}_m$ (for the
case of $f$ see \cite {Do}). Of course the main
difference between both functions is that $g^*$ is
intuitively total, while the whole point of the
$P=?NP$ question is whether $f$ is total or not. 

\subsection*{Summary of the paper}

$g^*$ is part of the folklore around the
$P=?NP$ question. As noticed, it is sometimes used as
a simile for the behavior of the counterexample
function to the $P=NP$ (no question--mark here!)
hypothesis. We consider its behavior in detail in
the present paper in order to show how wildly
different is the behavior of $g^*$ or its
modification $g$ (see below) when compared to that of
${\sf g}_0$. 

This paper shows that this function $g^*$ isn't
well--behaved at all. In fact, given a natural
version of it defined over all Turing machines, the
modified $g^*$ (noted $g$) will oscillate in a wild
way, and in its ``ups'' it will tower over all total
recursive functions, so that the only immediate
bounds we can find for it are in the realm of those
functions which bound all total recursive functions. 

One of those upper bounds is a generalized Busy
Beaver function, as we show in Section \ref {BB} of
this paper. 

\subsection*{Preliminary comments}

\begin {rem}\label {beg}\rm  If ${\cal M}$ is the set
of all Turing machines given by their G\"odel
numbers, the subset of all poly machines,
${\cal P}\subset{\cal M}$, isn't recursive. 

If we consider a function $g$ given by:
\begin {itemize}
\item $g(m) = \mu _x[{\sf M}_m(x) < 2^x]$, for
$m$ in ${\cal P}$, 
\item $g(m) = 0$, otherwise,  
\end {itemize} we intuitively see that $g$ is total.
\Vox
\end {rem}

We consider $g$ instead of $g^*$ as it is clear that
if there is some total recursive ${\sf f}$ that
bounds $g$, then we can obtain an adequately
modified ${\sf h}^*$ that will bound $g^*$. 

\subsection*{Goals}

Our goal in the present paper is to discuss some
properties of $g$. We will try to understand its
relation to the set of all total recursive
functions; everything proceeds within in an
intuitive framework, where we will be able to
compare it to a generalization of the Busy Beaver
function
\cite {Rado}. 

\

As noticed above, functions $g^*$ and $g$ are part
of the folklore around the $P=?NP$ question; $g^*$
was suggested to the authors by F.\ Cucker and,
independently, by W.\ Mitchell. 

\markboth {da Costa, Doria}{Busy Beaver}
\section{Conventions on Turing and
$\ell$--machines}\label {ppol}
\markboth {da Costa, Doria}{Busy Beaver}

We will use here what we call $\ell$--machines, for
``labeled'' or ``parametric'' Turing machines; they
are defined in this section. They are defined out of
the usual Turing machines, which we now describe.
They are simply Turing machines of which we keep
track with the help of a tag given by a parameter. 

\begin {rem}\rm\label {can} Suppose given the
canonical enumeration of binary words
$$\emptyset, 0, 1, 00, 01, 10, 11, 000, 001,\ldots$$
which code the empty word and the integers; they
correspond to 
$$0, 1, 2,\ldots.$$  We take this
correspondence to be fixed for the rest of this
paper. \Vox
\end {rem}

\subsection*{Turing machines}

\begin {rem}\label {TM}\rm We describe the behavior
of the Turing machines we deal here to avoid
ambiguities:
\begin {enumerate}
\item Turing machines are defined over the set
$A_2^*$ of finite words on the binary alphabet
$A_2=\{0,1\}$. 
\item Each machine has $n\geq 0$ states $s_0,
s_1,\ldots, s_{n-1}$, where $s_0$ is the final
state. (The machine stops when it moves to
$s_0$.)

We allow for a machine with $0$ states and an empty
table; it is discussed below in Remark
\ref {0}.  
\item The machine's head roams over a two--sided
infinite tape. 
\item Machines input a single binary word and either
never stop or stop when the tape has a finite, and
possibly empty set of binary words on it.  
\item The machine's {\it output word} will be the
one over which the head rests if and when
$s_0$ is reached. (If the head lies on a blank
square, then we take the output word to be the empty
word, that is, $0$.)  \Vox
\end {enumerate}
\end {rem}

\begin {rem}\rm \label {code} The Turing machine
inputs a binary string $\lfloor x\rfloor$ and (if it
stops over $\lfloor x\rfloor$) outputs a binary
string $\lfloor y\rfloor$. The corresponding
recursive function inputs the numeral $x$ and
outputs the numeral $y$. Whenever it is clear from
context, we write
$x$ for both the binary sequence and the numeral.
\Vox
\end {rem}

\begin {rem}\rm\label {ss} We will use upper case
and lower case sans serif letters (such as {\sf
M},\ldots) for Turing machines (and also for
$\ell$--machines). If
${\sf M}_n$ is a Turing machine of G\"odel number
$n$, its input--output relation is noted
${\sf M}_n(x)=y$. \Vox
\end {rem}

\begin {rem}\label{tour}\rm  Turing machines are
given by tables. We can write the tables as code
lines
$\xi,\xi',\ldots$, separated by blanks
$\sqcup$, such as 
$\xi\sqcup\xi'\sqcup\ldots\sqcup\xi''$. 

Let $\Xi$ be one such set of code lines separated by
blanks. Let $\Xi '$ be obtained out of $\Xi$ by a
permutation of the lines
$\xi,\xi ',\ldots$. Both $\Xi$ and $\Xi '$ are seen
as different machines that compute the same
algorithmic function, that is, in this case, if
$f_{\Xi}$ ($f_{\Xi '}$) is computed by $\Xi$ ($\Xi
'$), then $f_{\Xi} = f_{\Xi '}$.
\Vox
\end {rem}

\begin {rem}\label {0}\rm We define the {\it empty}
or {\it trivial machine} to be the Turing machine
with an empty table; we take it to be the simplest
example of the identity machine, again by
definition. (See also Remark \ref {0'}.) \Vox
\end {rem}

\subsection*{$\ell$--machines}

\begin {rem}\rm\label {ell}
We describe here the $\ell$--machines. Consider a
two--tape Turing machine \cite {Hop} where the input
is written over tape 1, while tape 2 comes with a
possibly empty binary string $n$. $n$ is the
machine's label, or parameter. 
\begin {itemize}
\item One easily sees that each Turing machine ${\sf
M}_m$ is simulated by infinitely many
$\ell$--machines: write an arbitrary $n$ on tape 2
and write all instructions for ${\sf M}_m$ solely
for tape 1. 

\item For the converse, if ${\sf M}_m$ is an
arbitrary Turing machine, and if (by an abuse of
notation) $\tau$ is the machine that polynomially
computes the onto pairing function $\omega \times
\omega\rightarrow\omega$, then the set of all
coupled machines ${\sf M}_m\circ \tau (n,x)$, all
$n, m$, represents the set of $\ell$--machines. 

\item Moreover, to avoid a pre--fixed label or
parameter, we can use the family of constant
Turing machines ${\sf i}_n$ that print $n$ over any
input, and couple it to the arrangement above, that
is ${\sf M}_m\circ \tau \circ {\sf i}_n$ (in order to
generate $n$), to obtain a set of emulated
$\ell$--machines. 
\end {itemize}

The G\"odel number $k = {\sf c}(m, \#\tau, n)$
($\#\ldots$ denotes the G\"odel number of $\ldots$)
is given by a primitive recursive ${\sf c}$ \cite
{MY}. \Vox
\end {rem}

>From that we have:

\begin {tef}
An $\ell$--machine is a pair $\langle n, {\sf
M}_m\rangle$, where $n\in\omega$ and ${\sf M}_m$ is
a 2--tape Turing machine of G\"odel number $m$, with
$n$ written on tape 2. \Vox
\end {tef}

\begin {rem}\rm
>From now on whenever we write ``machine'' we will
mean ``$\ell$--machine,'' unless otherwise stated.
\Vox
\end {rem}

\begin {rem}\rm
So, $\ell$--machines are just a convenient
bookkeeping device to keep track of families of
Turing machines which are labeled by some parameter.
However, under the guise of $\ell$--machines they
allow us the simple result stated in Lemma \ref
{koo}. \Vox
\end {rem}

\subsection*{G\"odel numbering}

(We use here Kleene's ``method of digits,'' without
however ``closing the gaps'' to squeeze out the
junk; see \cite {ML}, p.\ 178.)

\begin {tef}
Each $\ell$--machine is coded by the pair $\langle
n, m\rangle$, where both $n, m$ range over the whole
of $\omega$. \Vox
\end {tef}

\begin {rem}\rm

We refer to $\langle n, m\rangle$ as the {\em
G\"odel number of the $\ell$--machine,} or
{\em $\ell$--G\"odel number,} or {\em
$\ell$--index}. In order to compute some values for
$g$ as we do here one needs to know the
$\ell$--index for some poly machines which are given
as $\ell$--machines. \Vox

\end {rem}

\subsubsection*{A lemma on G\"odel numbers}

In order to compute pieces of $g$ one needs to
know the G\"odel numbers for some polynomial Turing
machines. We use the following result: 

\begin {lem}\label {koo} Let ${\sf M}_{\langle n,
m\rangle}$, $n\in\omega$, $m$ fixed, be a family of
$\ell$--machines. We can explicitly construct an
$\ell$--G\"odel numbering for the set of all
$\ell$--machines so that the $\ell$--indexes for
the ${\sf M}_{\langle n, m\rangle}$, $n\in\omega$,
are given by a linear function $N(n) = an + b$.  

Moreover there is a primitive recursive map from that
$\ell$--G\"odel numbering to standard G\"odel
numberings. Therefore, if ${\sf k}$ is such a map,
the induced $N'(n) = {\sf k}\circ N(n)$ is primitive
recursive. 
\end {lem} 

\begin {rem}\rm
Recall that the length $|x|$ of a finite string $x$
is the number of letters in $x$. \Vox
\end {rem}

{\it Proof of Lemma \ref {koo}}\,: The idea goes as
follows: $1547, 2547, 3547,
\ldots$ are in arithmetic progression with $r = 1000
= 10^3 = 10^{{\rm length}(547)}$. We extend this to
the proof of our lemma. 

\begin {itemize}

\item Recall the ordering of binary words as
described in Remark \ref {can}. 

\item $\ell$--machines are fully
characterized by pairs $\langle n, m\rangle$,
where $n$ is the parameter and $m$ the 2--tape 
Turing machine's code. 

\item If $x_n$ and $x_m$ are the corresponding
binary words (numerals), code pair
$\langle n, m\rangle$ as the word
$x_n10\widehat {x_m}$, where if $$x_m = x^0_m x^1_m
x^2_m\ldots x^k_m,$$ then
$$\widehat {x_m} = x^0_m x^0_m x^1_m x^1_m\ldots
x^k_m x^k_m.$$

\item So, given an arbitrary binary word, there is a
simple recursive procedure to decide whether it
codes an $\ell$--machine $\langle n,
m\rangle$ or not: go to the rightmost end and (by
going backwards) see if we have a duplicated word
that ends in the sequence $10$. If so, we have a
coded binary word for one of our pairs. 

\item Map the in--between mumbo--jumbo onto the
trivial $\ell$--machine.

\item This coding is certainly recursive and 1--1
onto the natural numbers. 

\item Finally, if $|x_m|$ is the length of $x_m$,
then
$|\widehat {x_m}| = 2|x_m| + 2$, and therefore,
between $\ell$--machine $\langle n, m\rangle$ and
$\ell$--machine
$\langle n + 1, m\rangle$ there are $2^{ 2|x_m| +
2}$ binary words. 
\end {itemize}

Follows the lemma, and the G\"odel number $N(n)$ of
the $\ell$--machines just described is given by an
arithmetic progression $N(n)$, of ratio $2^{ 2|x_m| +
2}$. 

>From Remark \ref {ell} one immediately sees that
there are primitive recursive maps that lead from
that coding to more usual G\"odel numberings, such
as the original one used by G\"odel, or Kleene's
``method of digits,'' so that $N'(n)$ is also
primitive recursive. \Vox

\begin {rem}\rm\label {Goed} So, it is enough to know
that the G\"odel numbers $N$ for the families of
machines we will be using here are related to the
machine's ``tag'' or parameter by $N = {\sf q}(n)$,
where
${\sf q}$ is primitive recursive, that is,
Ackermann's function ${\sf F}_{\omega}$ dominates
${\sf q}$. \Vox
\end {rem} 

\subsection*{Quasi--trivial machines}

We use $\ell$--machines; therefore here ``machine''
stands for $\ell$--machine. 

\begin {rem}\rm \label {0'}
Again the trivial $\ell$--machine is the one with the
empty table; more formally, $\ell$--machine $\langle
0, 0\rangle$. \Vox
\end {rem}

We will now consider a set of, say, nearly trivial
machines. They are all polynomial. Recall that the
operation time of a Turing  machine is given as
follows: if {\sf M} stops over an input
$x$, then the operation time over $x$, 
$$t_{\sf M} = |x| + \mbox {number of cycles of the 
machine until it stops.}$$ 

\begin {exa}\rm\ 

\begin {itemize}
\item {\bf First nearly trivial machine.} Note it
${\sf O}$.
${\sf O}$ inputs $x$ and stops. $$t_{\sf O} = |x| +
\mbox {moves to halting state} + \mbox {stops}.$$ 

So, operation time of ${\sf O}$ has a linear bound. 

\item {\bf Second nearly trivial machine.} Call it
${\sf O}'$. It  inputs $x$, always outputs $0$
(zero) and stops. 

Again operation time of ${\sf O}'$ has a linear
bound. 

\item {\bf Quasi--trivial machines.} A {\it
quasi--trivial machine} ${\sf Q}$ operates as
follows: for $x\leq  x_0$,
$x_0$ a constant value, ${\sf Q} = {\sf R}$,
${\sf R}$ an arbitrary total machine. For $x > x_0$,
${\sf Q} = {\sf O}$ or ${\sf O'}$. 

This machine has also a linear bound. \Vox
\end {itemize}
\end {exa}

We will use several $\ell$--families of
quasi--trivial machines. Please allow for some abuse
of language here; we need it in order to avoid 
cumbersome notations. (For instance, the machines
defined below depend on the G\"odel numbers of the
machines that appear as their subroutines.) 

\begin {rem}\rm\label {qt} Now let ${\sf H}$ be any
fast--growing, superexponential total 
machine. Let ${\sf H}'$ be another such machine. Form
the following family
${\sf Q}^{{\sf H}(n)}$ of quasi--trivial
$\ell$--machines with subroutines ${\sf H}$ and ${\sf
H}'$: 
\begin {enumerate}
\item ${\sf H}(n) = k(n)$, all $n$, is the way we
introduce the parameter in the family. 
\item If $x\leq k(n)$, ${\sf Q}^{{\sf H},{\sf
H'},n}(x) = {\sf H}'(x)$;
\item If $x > k(n)$, ${\sf Q}^{{\sf H},{\sf H'},n}(x)
= 0$.
\Vox
\end {enumerate}
\end {rem}

For function $g$ in Remark \ref {beg}: 

\begin {prop} If $N(n)$ is the G\"odel number of an 
$\ell$--machine as in Remark \ref {qt}, then $g(N(n))
= k(n) + 1 = {\sf H}(n) + 1$. 
\end {prop} 

{\it Proof}\,: Use Lemma \ref {koo}. Recall that
$N(n) = an + b$, $a$ and $b$ constants that depend
on the actual algorithms (i.e., in the G\"odel
numbers) for ${\sf H}, {\sf H}'$ and ${\sf Q}$. \Vox

\markboth {da Costa, Doria}{Busy Beaver}
\section {A domination lemma}
\markboth {da Costa, Doria}{Busy Beaver}

Recall:

\begin {tef}\label {dom} For
$f,g:\omega\rightarrow\omega$,
$$f \mbox{{\bf\ dominates\ }} g\leftrightarrow _{\rm
Def}
\exists y\,\forall x\,(x > y\rightarrow f(x)
\geq g(x)).$$ We write $f\succ g$ for $f$ dominates
$g$. \Vox
\end {tef} 

Our goal here is to prove the following result:

\begin {prop} For no total recursive function
${\sf h}$ does
${\sf h}\succ g$. 
\end {prop}

{\it Proof}\,: Suppose that there is a total
recursive function ${\sf h}$ such that
${\sf h}\succ g$. 

\begin {rem}\label {lev}\rm  Given such a function
${\sf h}$, obtain another total recursive function
${\sf h}'$ which satisfies:
\begin {enumerate}
\item ${\sf h}'$ is strictly increasing.
\item For $n > n_0$, ${\sf h}'(n) > {\sf h}(\frac
{a}{b}n -
\frac {1}{b})$, for integer values of the argument
of ${\sf g}$. 

That is, ${\sf h}'(an + b) > {\sf h}(n)$. \Vox
\end {enumerate}
\end {rem}

Constants $a$, $b$ are from the G\"odel numbers of
the quasi--trivial machines described in Remark
\ref {qt}. 

\begin {lem}\label {domi} Given a total recursive
${\sf h}$, there is a total recursive
${\sf h}'$ that satisfies the conditions in Remark
\ref {lev}. 
\end {lem} 

{\it Proof}\,: Given ${\sf h}$, obtain out of that
total recursive function a strictly increasing total
recursive
${\sf h}^*$. Then if, for instance, ${\sf
F}^{\omega}$ is Ackermann's function,
${\sf h}'={\sf h}^*\circ {\sf F}^{\omega}$ will do.
\Vox

\

>From Lemma \ref {koo}, we have that the G\"odel
numbers
$\#[{\sf Q}^{{\sf h}',{\sf K},n}]$ of the ${\sf
Q}^{{\sf h}',{\sf K},n}$ are given by $\#[{\sf
Q}^{{\sf h}',{\sf K},n}] = an + b$, $a, b\in\omega$. 

Therefore, $g(an + b) = {\sf h}'(n) + 1$. From Lemma
\ref {koo}, Remark \ref {lev} and Lemma \ref {domi}
we conclude our argument. If we make explicit the
computations, for
${\sf q}(n) = an+b$ (as the argument holds for any
strictly increasing primitive recursive ${\sf q}$): 
$$g({\sf q}(n)) = {\sf h}'(n) + 1 = {\sf h}^*({\sf
F}_{\omega}(n)) + 1,$$ and
$${\sf h}^*({\sf F}_{\omega}(n)) > {\sf h}^*({\sf
q}(n)).$$ For $N = {\sf q}(n)$ (see Remark \ref
{Goed}), 
$$g(N) > {\sf h}^*(N)\geq {\sf h}(N), \mbox {\
all $N$}.$$ 

Therefore no such ${\sf h}$ can dominate $g$. \Vox

\begin {rem}\rm
No need to emphasize that we could have used
ordinary Turing machines, as at the end of Remark
\ref {ell}; then $N(n)$ would be primitive recursive
on $n$. \Vox
\end {rem}

\markboth {da Costa, Doria}{Busy Beaver}
\section{An extended Busy Beaver function}\label {BB}
\markboth {da Costa, Doria}{Busy Beaver}

So $g$ overtakes all total recursive functions. Can
we put a ceiling to it? Yes: and that ceiling is a
kind of generalized Busy Beaver function (which
depends on $g$ itself). 

We now sketch an argument that shows that this
generalized Busy Beaver function dominates
$g$. 

The Busy Beaver function $B$ can be easily and
intuitively defined. If $N$ is the number of states
of a given Turing machine, then $B(N)$ is the
biggest number of $1$s a
$N$--state Turing machine prints over a tape filled
with $0$s as its input. 

\begin {rem}\rm  We are going to generalize that
function to
$B'$ as follows:
\begin {itemize}
\item First, notice that there is a simple relation
between a Turing machine's table and the number
\cite {Herm} of its states. 

More precisely, from that reference we see that a
$N$--state Turing machine with a binary alphabet has
a table given by a $N' = (3N + 1) \times 4$ matrix. 

So we can define a primitive recursive relation
${\sf f}$ from the G\"odel numbers to the set of
states of Turing machines so that if G\"odel numbers
$m'>m$, then states ${\sf f}(m')\geq {\sf f}(m)$. 

\item We can therefore see the Busy Beaver function
$B$ as a monotonic nondecreasing function of Turing
machine G\"odel numbers. 

\item Now for each G\"odel number $m$ for a {\it
polynomial} Turing machine, $g(m)$ is effectively
computed. 

\item We will define $B'(m)$ as: the maximum number
of $1$s printed by a $N(m)$--state Turing machine
over a tape with only $0$s and over a tape with all
inputs $< g(m) + 1$. 

\item Clearly $B'(N) \geq B(N)$ and $B'(N(m)) =
B(N(m))$ if the Turing machine with G\"odel number
$m$ isn't a polynomial machine. \Vox

\end {itemize}
\end {rem}

\begin {prop} For all $m$, $B'(m) \geq B(m)$. \Vox
\end {prop}

\begin {prop}
$B'$ dominates $g$. \Vox
\end {prop} 

\begin {rem}\rm We conjecture that $g$ is of the
same order of growth as the Busy Beaver function
itself. \Vox
\end {rem}

We intend to apply those ideas to the counterexample
function $f$ in a follow--up paper. 

\markboth {da Costa, Doria}{Busy Beaver}
\section {Acknowledgments}
\markboth {da Costa, Doria}{Busy Beaver}

The present paper is part of the Research Project on
Complexity and the Foundations of Computation, at
the Institute for Advanced Studies, University of
S\~ao Paulo (IEA--USP). 

The authors thank the Institute for support,
especially its Director Prof.\ A.\ Bosi, as well as
F.\ Katumi and S.\ Sedini. FAD also wishes to thank
Professors S.\ Amoedo de Barros and A.\ Cintra at
the Federal University of Rio de Janeiro, as well
its Rector Jos\'e Vilhena. 

Portions of the main questions dealt with in the
present paper were discussed and considered in the
newsgroup {\sc theory--edge@yahoogroups.com} from
November to December 2000; please search its files
for details. We heartily thank its participants for
criticisms and suggestions. 

\markboth {da Costa, Doria}{Busy Beaver}

\bibliographystyle {plain}

\end {document}